\documentclass[11pt]{article}%
\usepackage{amsmath}
\usepackage{amssymb}
\usepackage{amsfonts}
\usepackage{a4wide}

\newtheorem{example}{Example}

\begin{document}

\title{On the paper ``A study on concave optimization via canonical dual function''}

\author{C. Z\u{a}linescu\thanks{University
{}``Al.I.Cuza'' Ia\c{s}i, Faculty of Mathematics, Romania, email:
\texttt{zalinesc@uaic.ro}.}}
\date{}
\maketitle

\begin{abstract} In this short note we prove by a counter-example that
Theorem 3.2 in the paper ``A study on concave optimization via
canonical dual function'' by J. Zhu, S. Tao, D. Gao is false;
moreover, we give a very short proof for Theorem 3.1 in the same
paper.
\end{abstract}


\bigskip

In \cite{Zhu/Tao/Gao:09} one says: \textquotedblleft The primary
goal of this paper is to study the global minimizers for the
following concave optimization problem (primal problem $(P)$ in
short).

$(P)$ $\min$ $P(x)\quad$(1.1)

s.t. $x\in D,$

\noindent where

$D=\{x\in R^{n}\mid\left\Vert x\right\Vert \leq1\}$

\noindent and $P(x)$ is a smooth function in $R^{n}$ and is strictly concave
on the unit ball $D$, i.e. $\nabla^{2}P(x)<0$, on $D$.\textquotedblright

Even if it is not said what is meant by ``smooth function'', from
the context we think that $P$ is assumed to be a $C^2$ function on
$R^{n}$. One continues with \textquotedblleft Let's consider the
equation

$\left\{
\begin{array}
[c]{l}%
\nabla P(x)+\rho^{\ast}x=0,\quad x^{T}x=1,\\
\rho^{\ast}>0.
\end{array}
\right.  $ $\quad$(2.1)

\noindent Suppose there are only finitely many of root pairs for (2.1):

$0<\rho_{1}^{\ast}<\rho_{2}^{\ast}<\cdots<\rho_{l}^{\ast},$

\noindent associated with feasible points on the unit sphere:

$\widehat{x}_{1},\widehat{x}_{2},\ldots,\widehat{x}_{l},$

\noindent such that for each $i$,

$%
\begin{array}
[c]{l}%
\nabla P(\widehat{x}_{i})+\rho_{i}^{\ast}\widehat{x}_{i}=0,\quad\widehat
{x}_{i}^{T}\widehat{x}_{i}=1,\\
\rho_{i}^{\ast}>0.
\end{array}
$ $\quad$(2.2)\textquotedblright

Moreover, one says: \textquotedblleft In Section 3, two sufficient conditions
for determining a global minimizer are presented.\textquotedblright\

The results of \cite{Zhu/Tao/Gao:09} are the following.

\medskip
\textquotedblleft\textbf{Theorem 3.1.} If
$\nabla^{2}P(x)+\rho_{l}^{\ast}I>0$ on $\left\Vert x\right\Vert
\leq1$, then $\widehat{x}_{l}$ is a global minimizer of
(1.1).\textquotedblright

\medskip
\textquotedblleft\textbf{Theorem 3.2.} Suppose for $i=1,2,\ldots,l,$
$\det\left[  \nabla^{2}P(\widehat{x}_{i})+\rho_{i}^{\ast}I\right]
\neq0$ and $\frac{d^{2}P_{d}(\rho_{i}^{\ast})}{d\rho^{\ast2}}>0$.
Then $\widehat{x}_{l}$ is a global minimizer of
(1.1).\textquotedblright

\medskip
Related to these results we mention that Theorem 3.1 is (almost)
trivial and Theorem 3.2 is false even for $n=1.$

Indeed, because $\nabla^{2}P(x)+\rho_{l}^{\ast}I>0$ for $x\in D,$
there exists $r>1$ such that $\nabla^{2}P(x)+\rho_{l}^{\ast}I>0$ for
$x\in D_{r}:=\{u\in R^{n}\mid\left\Vert u\right\Vert <r\}.$
Otherwise there exist the sequences $(x_{k})\subset R^{n}$ with
$1<\left\Vert x_{k}\right\Vert \rightarrow1$ and $(v_{k})\subset
S=\{u\in R^{n}\mid\left\Vert u\right\Vert =1\}$ such that
$v_{k}^{T}\left(  \nabla^{2}P(x_{k})+\rho_{l}^{\ast}I\right)
v_{k}\leq0$ for every $k$. We may assume that $x_{k}\rightarrow x$
and $v_{k}\rightarrow v;$ hence $x,v\in S$. It follows that
$v^{T}\left(  \nabla^{2}P(x)+\rho_{l}^{\ast }I\right)  v\leq0,$
contradicting our assumption. Since $D_{r}$ is an open convex set we
obtain that $P+\tfrac{1}{2}\rho_{l}^{\ast}\left\Vert
\cdot\right\Vert ^{2}$ is a (strictly) convex function on $D_{r}.$
Because $\widehat{x}_{l}\in D\subset D_{r}$ and $\nabla\big(
P+\tfrac{1}{2}\rho_{l}^{\ast}\left\Vert \cdot\right\Vert ^{2}\big)
(\widehat{x}_{l})=\nabla
P(\widehat{x}_{l})+\rho_{l}^{\ast}\widehat{x}_{l}=0,$
we have that $\widehat{x}_{l}$ is a global minimizer of $P+\tfrac{1}{2}%
\rho_{l}^{\ast}\left\Vert \cdot\right\Vert ^{2}$ on $D_{r}$. In
particular we have that
\[
P(\widehat{x}_{l})+\tfrac{1}{2}\rho_{l}^{\ast}=P(\widehat{x}_{l})+\tfrac{1}%
{2}\rho_{l}^{\ast}\left\Vert \widehat{x}_{l}\right\Vert ^{2}\leq
P(x)+\tfrac{1}{2}\rho_{l}^{\ast}\left\Vert x\right\Vert ^{2}\leq
P(x)+\tfrac{1}{2}\rho_{l}^{\ast}\quad\forall x\in D,
\]
whence $P(\widehat{x}_{l})\leq P(x)$ for every $x\in D.$

The proof above shows that whenever $P$ is a $C^{2}$ function on an open set
$D_{r}$ containing $D$ such that $\nabla^{2}P(x)+\rho_{l}^{\ast}I>0$ on $D$
(or even less, $\nabla^{2}P(x)+\rho_{l}^{\ast}I\geq0$ on $D_{r}$) and
$\overline{x}\in S$ and $\overline{\rho}\geq0$ are such that $\nabla
P(\overline{x})+\overline{\rho}\,\overline{x}=0$, then $\overline{x}$ is a
global minimizer of $P$ on $D.$

\bigskip

Related to \cite[Th.\ 3.2]{Zhu/Tao/Gao:09}, let us observe first
that the condition
$\frac{d^{2}P_{d}(\rho_{i}^{\ast})}{d\rho^{\ast2}}>0$ is equivalent
with $\widehat{x}_{i}^{T}\left[
\nabla^{2}P(\widehat{x}_{i})+\rho_{i}^{\ast }I\right]
^{-1}\widehat{x}_{i}<0.$

Indeed, one says: \textquotedblleft For $i=1,2,\ldots,l$, defined by

$\nabla P(\widehat{x}(\rho^{\ast}))+\rho^{\ast}\widehat{x}(\rho^{\ast
})=0,\quad\rho^{\ast}>0,\ \ \widehat{x}(\rho_{i}^{\ast})=\widehat{x}_{i}\quad$ (2.3)

\noindent a branch $\widehat{x}_{i}(\rho^{\ast})$ is a continuously
differentiable vector function on $\rho^{\ast}$.\textquotedblright%
\ \textquotedblleft In what follows, we suppress the index when focusing on a
given branch according to the context.

The dual function [6] with respect to a given branch $\widehat{x}(\rho^{\ast
})$ is defined as

$P_{d}(\rho^{\ast})=P(\widehat{x}(\rho^{\ast}))+\frac{\rho^{\ast}}{2}%
\widehat{x}^{T}(\rho^{\ast})\widehat{x}(\rho^{\ast})-\frac{\rho^{\ast}}%
{2}.\quad$ (2.6)\textquotedblright

\smallskip
Note that [6] above is our reference \cite{Gao-book}.

\smallskip
In order to obtain a solution $\widehat{x}$ of (2.3) the authors use
differential equations. In fact, let $F:R^{n}\times R\rightarrow
R^{n}$ be defined by $F(x,\rho):=\nabla P(x)+\rho x$. Clearly, $F$
is a $C^{1}$ function, $\nabla_{\rho}F(x,\rho)=x$, whence
$\nabla_{\rho}F(\widehat{x}_{i},\rho_{i}^{\ast})=\widehat{x}_{i}\neq0.$
By the implicit function theorem a $C^{1}$ function
$\widehat{x}:J\rightarrow R^{n}$ exists such that
$F(\widehat{x}(\rho),\rho)=0$ for $\rho\in J$ and
$\widehat{x}(\rho_{i}^{\ast})=\widehat{x}_{i}$, where $J$ is an open
interval containing $\rho_{i}^{\ast}$. It follows that
\[
\nabla_{x}F(\widehat{x}(\rho),\rho)\widehat{x}^{\prime}(\rho)+\nabla_{\rho
}F(\widehat{x}(\rho),\rho)=\left[  \nabla_{x}^{2}P(\widehat{x}(\rho))+\rho
I\right]  \widehat{x}^{\prime}(\rho)+\widehat{x}(\rho)=0\quad\forall\rho\in J.
\]
Because $\det\left[
\nabla^{2}P(\widehat{x}_{i})+\rho_{i}^{\ast}I\right] \neq0$, we may
assume that $\det\left[  \nabla^{2}P(\widehat{x}(\rho
))+\rho_{i}^{\ast}I\right]  \neq0$ for all $\rho\in J$ (taking a
smaller $J$
if necessary). Hence%
\[
\widehat{x}^{\prime}(\rho)=-\left[  \nabla_{x}^{2}P(\widehat{x}(\rho))+\rho
I\right]  ^{-1}\widehat{x}(\rho)\quad\forall\rho\in J.
\]
From the expression of $P_{d}$ in (2.6), using (2.3) we get
\begin{align*}
P_{d}^{\prime}(\rho) &  =\nabla P(\widehat{x}(\rho))\widehat{x}^{\prime}%
(\rho)+\tfrac{1}{2}\widehat{x}^{T}(\rho)\widehat{x}(\rho)+\rho\widehat{x}%
^{T}(\rho)\widehat{x}^{\prime}(\rho)-\tfrac{1}{2}=\tfrac{1}{2}\widehat{x}%
^{T}(\rho)\widehat{x}(\rho)-\tfrac{1}{2},\\
P_{d}^{\prime\prime}(\rho) &  =\widehat{x}^{T}(\rho)\widehat{x}^{\prime}(\rho)
\end{align*}
for every $\rho\in J$. Using the expression of
$\widehat{x}^{\prime}(\rho)$ obtained above we get
\[
P_{d}^{\prime\prime}(\rho_{i}^{\ast})=-\widehat{x}^{T}(\rho_{i}^{\ast})\left[
\nabla_{x}^{2}P(\widehat{x}(\rho_{i}^{\ast}))+\rho_{i}^{\ast}I\right]
^{-1}\widehat{x}(\rho_{i}^{\ast})=-\widehat{x}_{i}^{T}\left[  \nabla
^{2}P(\widehat{x}_{i})+\rho_{i}^{\ast}I\right]  ^{-1}\widehat{x}_{i}.
\]

This shows that instead of the condition $\frac{d^{2}P_{d}(\rho_{i}^{\ast}%
)}{d\rho^{\ast2}}>0$, which uses a quite complicated function, it
was preferable to consider the condition $$\widehat{x}_{i}^{T}\left[
\nabla^{2}P(\widehat {x}_{i})+\rho_{i}^{\ast}I\right]
^{-1}\widehat{x}_{i}<0,$$ which is written using the data of the
problem.

\begin{example}
Consider $P:R\rightarrow R$ defined by
$P(x)=-x^{4}-\frac{8}{5}x^{3}-\frac {6}{5}x^{2}+\frac{12}{5}x$. Then
\[
\textstyle P^{\prime}(x)=-4x^{3}-\frac{24}{5}x^{2}-\frac{12}{5}x+\frac{12}%
{5},\quad P^{\prime\prime}(x)=-12x^{2}-\frac{48}{5}x-\frac{12}{5}.
\]
We have that $P^{\prime\prime}(x)\leq P^{\prime\prime}(-\frac{2}{5}%
)=-\frac{12}{25}<0$ for every $x\in R;$ hence $P$ is a strictly
concave function. The system (2.1) becomes $x=\pm1,$
$\rho^{\ast}=-x^{-1}P^{\prime }(x)$, $\rho^{\ast}>0$. The solutions
are $(\widehat{x}_{i},\rho_{i}^{\ast})$, $i\in\{1,2\}$, where
$\widehat{x}_{1}=-1$, $\widehat{x}_{2}=1$, $\rho
_{1}^{\ast}=P^{\prime}(-1)=4,$
$\rho_{2}^{\ast}=-P^{\prime}(1)=\frac{44}{5}$. Hence $l=2$ and
$0<\rho_{1}^{\ast}<\rho_{2}^{\ast}$. The condition
$\widehat{x}_{i}^{T}\left[
\nabla^{2}P(\widehat{x}_{i})+\rho_{i}^{\ast
}I\right]  ^{-1}\widehat{x}_{i}<0$ becomes $P^{\prime\prime}(\widehat{x}%
_{i})+\rho_{i}^{\ast}<0$ for $i\in\{1,2\}$, in which case
$\det\left[ \nabla^{2}P(\widehat{x}_{i})+\rho_{i}^{\ast}I\right]
\neq0$. But
$P^{\prime\prime}(-1)+4=-\frac{4}{5}<0$, $P^{\prime\prime}(1)+\frac{44}%
{5}=-\frac{76}{5}<0$. Using \cite[Th.\ 3.2]{Zhu/Tao/Gao:09} we
obtain that $\widehat{x}_{2}=1$ is the global minimizer of $P$ on
$[-1,1].$ However, $P(-1)=-3<-\frac{7}{5}=P(1)$, proving that
\cite[Th.\ 3.2]{Zhu/Tao/Gao:09} is false.
\end{example}

\end{document}